\documentclass[]{amsart}
\usepackage{amssymb}
\usepackage{graphicx}
\usepackage{xcolor} 
\usepackage{amsmath}
\usepackage{amsthm}
\usepackage{verbatim}
\usepackage{enumitem}

\newtheorem{theorem}{Theorem}

\newtheorem{corollary}[theorem]{Corollary}
\newtheorem*{theorem*}{Theorem}
\newtheorem*{lemma*}{Lemma}

\theoremstyle{definition}
\newtheorem{definition}{Definition}
\newtheorem{remark}{Remark}

\newtheorem*{remark*}{Remark}
\newtheorem*{example*}{Example}
\newtheorem*{er*}{Examples and Remarks}

\newcommand{\nnrm}[1]{{\vert\kern-0.25ex\vert\kern-0.25ex\vert #1 \vert\kern-0.25ex\vert\kern-0.25ex\vert}}

\begin{document}

\title[Locality of vortex stretching]
{Locality of vortex stretching for the 3D Euler equations} 

\author{Yuuki Shimizu} 
\address{Graduate School of Mathematical Sciences, University of Tokyo, Komaba 3-8-1 Meguro, 
Tokyo 153-8914, Japan } 
\email{yshimizu@g.ecc.u-tokyo.ac.jp} 

\author{Tsuyoshi Yoneda} 
\address{Graduate School of Economics, Hitotsubashi University, 2-1 Naka, Kunitachi, Tokyo 186-8601, Japan} 
\email{t.yoneda@r.hit-u.ac.jp}


\date{\today} 


\begin{abstract} 
We consider 
the 3D incompressible Euler equations
 under the following situation: small-scale  vortex blob being stretched by a prescribed large-scale stationary flow. 
More precisely, we clarify what kind of large-scale stationary flows really stretch small-scale vortex blobs in alignment with the straining direction.
The key idea is constructing a Lagrangian coordinate so that the Lie bracket is identically zero (c.f. the Frobenius theorem), and investigate the locality of  the pressure term by using it.
\end{abstract}

\maketitle

\section{Introduction}

The most important features of the Navier-Stokes turbulence is that turbulence is not random but composed of vortex stretching behavior. More precisely, 
recent DNS
\cite{Goto-2008,GSK,Motoori-2019,Motoori-2021} of
the Navier-Stokes turbulence at sufficiently high Reynolds numbers have reported that
there exists a hierarchy of vortex stretching motions.
In particular, Goto-Saito-Kawahara \cite{GSK} clearly observed that turbulence at sufficiently high Reynolds numbers in a periodic cube is composed of a self-similar
hierarchy of antiparallel pairs of vortex tubes, and  it is sustained
by creation of smaller-scale vortices due to stretching in larger-scale strain fields. 
This observation is further investigated by Y-Goto-Tsuruhashi \cite{YGT} (see also \cite{TGOY}). Thus we could conclude physically  
that  local-scale energy transfer is mainly induced by vortex stretching (see also \cite{JY1,JY2,JY3} for the related mathematical results).
Therefore as the sequence of these studies, our next study will be clarifying the vortex stretching dynamics precisely.

In this paper we mathematically consider the locality of small-scale vortex dynamics in the 3D incompressible Euler equations. 
More precisely, 
we consider the inviscid flow under the following situation: small-scale  vortex blob being stretched by a prescribed
large-scale stationary flow, and 
we clarify what kind of  large-scale stationary flows really stretch
 smaller-scale vortex blobs in alignment with the straining direction. 
Now let us describe the incompressible Euler equations (inviscid flow)  as follows:
\begin{eqnarray}
\label{Euler eq.}
& &
\partial_tu+(u\cdot \nabla)u=\left(\partial_t(u\circ \Phi)\right)\circ\Phi^{-1}=
\partial_t^2\Phi\circ\Phi^{-1}=-\nabla p,\quad\nabla \cdot u=0\quad \text{in}
\quad \mathbb{R}^3,\\
\nonumber
& & \quad u|_{t=0}=u_0,
\end{eqnarray}
where $\Phi$ is 
 the associated Lagrangian flow given by \begin{equation*}
\begin{split}
\partial_t \Phi(t,x)=u(t,\Phi(t,x))=:u\circ \Phi\quad\text{with}\quad \Phi(0,x)=x\in\mathbb{R}^3.
\end{split}
\end{equation*} 
Let $u^S:(-\epsilon,\epsilon)\times\mathbb{R}^3\to\mathbb{R}^3$  be the flow of a small-scale vortex blob and $u^L:\mathbb{R}^3\to\mathbb{R}^3$ be a prescribed large-scale stationary flow. Then
the associated Lagrangian flows $\eta^L$ and $\eta^S$ satisfying $\Phi(t,x)=\eta^L(t,\eta^S(t,x))=\eta^L\circ\eta^S$ are given by 
\begin{equation*}
\begin{split}
&\partial_t \eta^L(t,x)=u^L(\eta^L(t,x))=:u^L\circ \eta^L\quad\text{with}\quad \eta^L(0,x)=x\in\mathbb{R}^3,\\
&\partial_t \eta^S(t,x)=u^S(t,\eta^S(t,x))=:u^S\circ \eta^S\quad\text{with}\quad \eta^S(0,x)=x\in\mathbb{R}^3.
\end{split}
\end{equation*}
Let us assume $\eta^S(t,x_0)=\eta^S(0,x_0)=x_0\in\ell$, where
\begin{equation}\label{alignment}
\ell:=\bigcup_{t\in(-\epsilon,\epsilon)}\eta^L(t,x_*)\quad\text{for some}\quad x_*\in\mathbb{R}^3.
\end{equation} 
$\ell$ represents the rotating axis of the small-scale vortex blob which  
 aligns
with the large-scale straining flow.
Based on DNS of
homogeneous isotropic turbulence,  Hamlington-Schumacher-Dahm \cite{HSD} showed vorticity tends to align with the stretching
direction of the background strain.  Note that,  
in their study, 
the background strain means the strain rate induced by the vorticity beyond radius $\sim12\eta$ ($\eta$ is the Kolmogorov
scale).
 See also \cite{GSK}.
With the aid of this physical observation, we constructed 
the mathematical model \eqref{alignment}.

In general, $\eta^S$ is strongly affected by the large-scale straining flow $\eta^L$ when $\Phi=\eta^L\circ\eta^S$ comes from an Euler flow, thus, in this study, 
we need to clarify the nonlinear interaction even partially. 
\begin{remark}\label{typical}
For readers' convenience, in this remark, we state a typical vortex stretching
as an example. Let $x_0\in\{x_1=x_2=0\}=\ell$, 
and let $(r,\theta,x_3)$ be the cylindrical coordinate such that $(x_1,x_2)=(r\cos\theta, r\sin\theta)$ with
\begin{equation*}
e_r:=\frac{(x_1,x_2,0)}{\sqrt{x_1^2+x_2^2}},\quad
e_\theta:=\frac{(-x_2,x_1,0)}{\sqrt{x_1^2+x_2^2}}\quad\text{and}\quad
e_{x_3}:=(0,0,1).
\end{equation*}
First we give the typical straining flow $(u_r^L, u^L_\theta, u_{x_3}^L)$ as follows:
\begin{equation*}
u^L_r:=u^L\cdot e_r=-r,\quad u^L_{x_3}:=u^L\cdot e_{x_3}=2x_3\quad\text{and}\quad u^L_\theta:=u^L\cdot e_\theta=0. 
\end{equation*}
In this case the corresponding Lagrangian flow is 
\begin{equation*}
(\eta^L_r,\eta^L_\theta, \eta^L_{x_3})=(e^{-t} r_0,\theta_0, e^{2t}x_{3,0}),
\end{equation*}
 where $(r_0, \theta_0,x_{3,0})$ represents
the initial position.
For $\omega:=\nabla\times u$,
let $\omega^S:=\nabla\times u^S$ be axisymmetric vorticity depending only on $r$ variable, such that 
\begin{equation*}
\omega^S_\theta(t,r):=\omega^S\cdot e_\theta,\quad
\omega^S_{x_3}:=\omega^S\cdot e_{x_3}=0\quad\text{and}\quad
\omega^S_r:=\omega^S\cdot e_r=0.
\end{equation*}
Then we have the following explicit solution (vorticity) to the Euler equations:
\begin{equation}\label{typical solution}
\omega_\theta^S(t,\eta^L_r(t))=e^t\omega^S_{0,\theta}(e^{-t}r_0)\quad\text{and}\quad \omega^S_{0,\theta}=\omega^S_\theta(t)|_{t=0}.
\end{equation}
\end{remark}
\vspace{0.3cm}

With the aid of the typical stretching motion in Remark \ref{typical},
we rigorously define the meaning of  ``stably stretching" as follows.
\begin{definition}(Stably stretching.)\quad 
 $\eta^L$
 is said to be \textit{stably stretching} along the rotating axis $\ell\subset \mathbb{R}^3$, if 
the following two conditions hold.
\begin{itemize}
\item The rotation axis $\ell$ is stretched along $\partial_t\eta^L$ direction, that is 
\begin{equation}\label{stretch}
\partial_t|\partial_t\eta^L(t,x_0)|>0\quad\text{for each}\quad x_0\in\ell.
\end{equation} 
\item
Let 
\begin{equation*}
x_0^\perp:=\{x\in\mathbb{R}^3: (x_0-x)\cdot \partial_t\eta^L(0,x_0)=0,\ 
|x_0-x|< \delta\}\quad\text{for each}\quad x_0\in\ell.
\end{equation*}
The time evolution of the surface $x_0^\perp$ (accompanied by the fluid particles) is always perpendicular to the stretching direction at $\eta^L(t,x_0)$  ($t>0$), that is, 
\begin{equation}\label{perp}
\partial_t\eta^L(t,x_0)\perp\Phi(t,x_0^\perp)\quad\text{for each}\quad x_0\in\ell.
\end{equation}
\end{itemize}
\end{definition}
\begin{remark}
 \eqref{typical solution} satisfies 
\eqref{stretch}
and 
\eqref{perp}.
\end{remark}

To state our theorem, we need to prepare ``curvature". 
First let us choose a point $x_0\in\ell$ (reference point) and fix it.
Identifying  $\frac{d}{dz}$ with $\partial_z$, we define
 $t(z)$ as
\begin{equation}\label{ODE}
\partial_zt>0,\quad |\partial_z\eta^L|:=|\partial_z\eta^L(t(z),x_0)|=1\quad\text{and}\quad
t(0)=0.
\end{equation}
In this case we immediately have  $\partial_z t=|\partial_t\eta^L|^{-1}$
and its inverse is $\partial_tz=|\partial_t\eta^L|$ (with the variables omitted).
Then define the unit tangent vector $\tau$ as 
\begin{equation*}
\tau(z)=\partial_z\eta^L(t(z),x_0),
\end{equation*}
the unit curvature vector $n$ as $\kappa n=\partial_z \tau$ with a curvature function $\kappa(z)>0$,
the unit torsion vector $b$ 
as : $b(z):=\tau(z)\times n(z)$, where $\times$ is the exterior product.
Without loss of generality, we can assume the torsion function $T(z)$ is positive
by choosing the orientation of the torsion vector $b(z)$.
Then we now state our main theorem.
\begin{theorem}\label{criterion}
Let $(\partial_t \Phi)\circ\Phi^{-1}=\partial_t(\eta^L\circ\eta^S)\circ(\eta^L\circ\eta^S)^{-1}$ be a solution to the Euler equations.
If $\eta^L$
is stably stretching along $\ell\subset\mathbb{R}^3$,
 then we have 
\begin{equation}\label{vortex-stretching-condition}
\partial_z(\kappa|\partial_t\eta^L|^2)
=0\quad\text{for each}\quad x_0\in\ell,
\end{equation}
in particular, we do not need any information of $\eta^S$,
\end{theorem}
\begin{remark}
Since $\partial_t|\partial_t\eta^L|>0$, 
we obtain a necessary condition $\partial_z\kappa(z)\leq 0$.
\end{remark}
Finally, we examine a typical straining flow $(-x_1,x_2,0)$ whether or not it satisfies the formula \eqref{vortex-stretching-condition}.
Note that we could also examine the straining flow $(-x_1,-x_2,2x_3)$ which is already appeared in Remark \ref{typical}, but in this case the formula 
becomes much more complicated, thus we leave it as a reader's exercise. 
\begin{corollary}\label{typical-strain-tensor}
If $u^L$ is a straining flow such that $u^L=(-x_1,x_2,0)$, then we have 
\begin{equation*}
\partial_z(\kappa |\partial_t\eta^L|^2)=-\frac{\tanh 2t}{\cosh 2t}
\quad\text{with}\quad
t=\frac{1}{2}\log\left(\frac{x_1}{x_2}\right).
\end{equation*}
This means that, 
if $\ell\not\subset\{x_1=0\}$ and is in the stretching region: $\ell\subset 
\{x: |x_2|>|x_1|\}$, then the pair of  $\eta^L$ and $\Phi=\eta^L\circ\eta^S$ does not satisfy \eqref{perp},
which implies that  $\eta^L$ is  
``unstably" stretching along $\ell$.
\end{corollary}

\section{Proof of Theorem \ref{criterion}}

First we rephrase the initial flat plane $x_0^\perp$ such that 
\begin{equation}\label{rephrase}
x_0^\perp=\{x_0+r_1 n(0)+ r_2 b(0): \sqrt{r_1^2+r_2^2}<\delta\}.
\end{equation}
For any initial particle on the plane $x=x_0+ r_1n(0)+ r_2 b(0)\in x_0^\perp$, 
$\Phi(t,x)$
is uniquely expressed as (we omit the change of variables) 
\begin{equation}\label{Phi}
\begin{split}
\Phi(t,x)
&=:\Phi(z,r_1,r_2)\\
&=\eta^L(t(z),x_0)
+Z(z,r_1,r_2)\tau(z)
+R_1(z,r_1,r_2)n(z)+R_2(z,r_1,r_2)b(z),
\end{split}
\end{equation} 
for sufficiently small $r_1$ and $r_2$, with  $Z(z,0,0)=0$, $Z(0,r_1,r_2)=0$, $R_1(0,r_1,r_2)=r_1$ and $R_2(0,r_1,r_2)=r_2$.
\begin{remark}
We can rephrase \eqref{perp} as follows:
\begin{equation}\label{perpendicular}
\begin{split}
\partial_{r_1}Z|_{r_1=r_2=0}=\partial_{r_2}Z|_{r_1=r_2=0}=0\quad\text{for}\quad z>0.
\end{split}
\end{equation}
This is due to the fact that 
\begin{equation*}
\Phi(t,x_0^\perp)=\bigcup_{\sqrt{r_1^2+r_2^2}<\delta}\Phi(z(t),r_1,r_2)
\end{equation*}
by \eqref{rephrase} and \eqref{Phi}, and 
\begin{equation*}
\partial_t\eta^L(t,x_0)\cdot \partial_{r_1}\Phi(z(t),r_1,r_2)|_{r_1=r_2}=
\partial_t\eta^L(t,x_0)\cdot \partial_{r_2}\Phi(z(t),r_1,r_2)|_{r_1=r_2}=0.
\end{equation*}
\end{remark}
Since 
 the corresponding Jacobian $\frac{\partial(R_1,R_2)}{\partial(r_1,r_2)}$
is clearly nonzero for sufficiently small $z>0$, so, by the inverse function theorem (for each $z>0$),
 we can rewrite the equality \eqref{Phi} as follows:
Let $R_1$ and $R_2$ be variables, and $r_1$ and $r_2$ be the corresponding inverse functions. For any particle $x=x_0+r_1(z,R_1,R_2)n(0)+r_2(z,R_1,R_2)b(0)\in x_0^\perp$,
\begin{equation}\label{explicit-formula} 
\begin{split}
\Phi(t,x)
&=:\Phi(z,R_1,R_2)\\
&=\eta^L(t(z),x_0)
+Z(z,r_1(z,R_1,R_2),r_2(z,R_1,R_2))\tau(z)
+R_1n(z)+R_2b(z)
\end{split}
\end{equation} 
for sufficiently small $|R_1|, |R_2|, z>0$.
Let us recall the Frenet-Serret formulas:
\begin{equation*}
\frac{d}{dz}
\begin{pmatrix}
\tau\\
n\\
b
\end{pmatrix}
=
\begin{pmatrix}
0&\kappa&0\\
-\kappa& 0& T\\
0& -T& 0
\end{pmatrix}
\begin{pmatrix}
\tau\\
n\\
b
\end{pmatrix}
.
\end{equation*}
By combining the Frenet-Serret formulas, \eqref{perpendicular} and \eqref{explicit-formula}, we have
\begin{equation}\label{Frenet-Serret}
\begin{cases}
\partial_{z}\Phi&=
\tau+R_1(Tb-\kappa \tau)-R_2 T n+\mathcal O(R_1^2,R_2^2),\\
\partial_{R_1}\Phi&= n+\mathcal O(R_1,R_2),\\
\partial_{R_2}\Phi&=b+\mathcal O(R_1,R_2),
\end{cases}
\end{equation}
where $\mathcal O$ is Landau's notation.
A direct calculation yields
\begin{equation*}
\begin{split}
\partial_z^2\Phi=
&\kappa n-R_1(T^2+\kappa^2 )n+R_1((\partial_zT)b-(\partial_z\kappa) \tau)\\
&\ -R_2T(-\kappa\tau+Tb)
-R_2(\partial_zT) n+\mathcal O(R_1^2,R_2^2)
\end{split}
\end{equation*}
and then
\begin{equation*}
\begin{split}
&\partial_{R_1}\partial_z\Phi|_{R_1=R_2=0}=Tb-\kappa \tau,\\
&\partial_{R_1}\partial_z^2\Phi|_{R_1=R_2=0}=-(T^2+\kappa^2)n+((\partial_zT)b-(\partial_z\kappa)\tau),\\
&\partial_{R_2}\partial_z\Phi|_{R_1=R_2=0}=-Tn,\\
&\partial_{R_2}\partial_z^2\Phi|_{R_1=R_2=0}=-T(-\kappa \tau+Tb)-(\partial_zT)n.
\end{split}
\end{equation*}
Thus
\begin{equation*}
(\partial_{R_1}\partial_z^2\Phi)\cdot\tau|_{R_1,R_2=0}=-\partial_z\kappa\quad\text{and}\quad
(\partial_{R_1}\partial_z\Phi)\cdot\tau|_{R_1,R_2=0}=-\kappa.
\end{equation*} 
On the other hand, by the Leibniz rule, we see
\begin{equation*}
\partial_t^2\Phi=\partial_z^2\Phi(\partial_tz)^2+\partial_z\Phi\partial_t^2z
\quad\text{with the variables omitted.}
\end{equation*}
Combining the facts $\partial_t^2z=\partial_t|\partial_t\eta^L|$ and $\partial_tz=|\partial_t\eta^L|$, we have 
\begin{eqnarray}\label{pressure-estimate}
& &
-\partial_{R_1}(\nabla p\cdot\tau)=
\partial_{R_1}(\partial_t^2\Phi\cdot\tau)=(\partial_{R_1}\partial_t^2\Phi)\cdot\tau=
-\kappa\partial_t|\partial_t\eta^L|
-\partial_z\kappa |\partial_t\eta^L|^2,\\
\label{b}
& &
-\partial_{R_2}(\nabla p\cdot\tau)=
\partial_{R_2}(\partial_t^2\Phi\cdot\tau)=(\partial_{R_2}\partial_t^2\Phi)\cdot\tau=
+T\kappa |\partial_t\eta^L|^2
\end{eqnarray}
for $R_1=R_2=0$.
Next we derive other formulae by using the Euler equations.
By the Leibniz rule, we see
\begin{equation*}
\kappa n=\partial_z^2\eta^L(t(z),x_0)=\partial_z(\partial_t\eta^L\partial_zt)
=\partial_t^2\eta^L(\partial_zt)^2+\partial_t\eta^L\partial_z^2t.
\end{equation*}
Combining 
$\partial_z^2t=\partial_z|\partial_t\eta^L|^{-1}=-|\partial_t\eta^L|^{-2}\partial_z|\partial_t\eta^L|
=-|\partial_t\eta^L|^{-3}\partial_t|\partial_t\eta^L|$,
we have 
\begin{equation*}
\partial_t^2\eta^L=|\partial_t\eta^L|^{2}\kappa n+\partial_t|\partial_t\eta^L|\tau.
\end{equation*}
By using the Euler equations, we have
\begin{eqnarray*}
-\nabla p\cdot\tau
&=&
\partial_t^2\Phi\cdot \tau
=
\partial_{t}|\partial_t\eta^L|,\\
-\nabla p \cdot n
&=&
\partial_t^2\Phi\cdot n=
\kappa|\partial_t\eta^L|^2,\\
-\partial_{z}(\nabla p\cdot n)
&=&
\partial_z\kappa|\partial_t\eta^L|^2+2\kappa\partial_t|\partial_t\eta^L|,\\ 
-\nabla p\cdot b
&=&0
\end{eqnarray*}
for $R_1=R_2=0$ with the change of variables $\circ\eta^L$ omitted again. On the other hand, from \eqref{Frenet-Serret},
 we have  the following inverse matrix: 
\begin{equation*}
\begin{pmatrix}
\tau\\
n\\
b\\
\end{pmatrix}
=
\begin{pmatrix}
(1-\kappa R_1)^{-1}& R_2T (1-\kappa R_1)^{-1}& -R_1 T(1-\kappa R_1)^{-1}\\
0& 1& 0\\
0& 0& 1
\end{pmatrix}
\begin{pmatrix}
\partial_{z}\Phi\\
\partial_{R_1}\Phi\\
\partial_{R_2}\Phi
\end{pmatrix}
\end{equation*}
with the higher order terms omitted since we finally take $R_1,R_2\to 0$.
Then we see 
\begin{equation*}
\begin{split}
\nabla p\cdot \tau
=&
(1-\kappa R_1)^{-1}(\nabla p\cdot \partial_z\Phi)\\
&+
R_2T (1-\kappa R_1)^{-1}(\nabla p\cdot\partial_{R_1}\Phi)-R_1 T(1-\kappa R_1)^{-1}
(\nabla p\cdot \partial_{R_2}\Phi)\\
=&
(1-\kappa R_1)^{-1} \partial_z(p\circ\Phi)\\
&+
R_2T (1-\kappa R_1)^{-1}\partial_{R_1} (p\circ\Phi)-R_1 T(1-\kappa R_1)^{-1}
\partial_{R_2}(p\circ\Phi).
\end{split}
\end{equation*}
and then (omit the variable $\circ\Phi$)
\begin{eqnarray*}
-\partial_{R_1}(\nabla p\cdot \tau)|_{R_1=R_2=0}
&=&
\left(-\kappa\partial_{z} p-\partial_{R_1}\partial_{z} p
-T\partial_{R_2}p\right)|_{R_1=R_2=0}\\
\nonumber
(\text{commute}\ \partial_{R_1}\ \text{and}\  \partial_{z})
&=&
\left(-\kappa(\nabla p\cdot \tau)-\partial_{z} (\nabla p\cdot n)-T(\nabla p\cdot b)\right)|_{R_1=R_2=0}\\
&=&
3\kappa \partial_t|\partial_t\eta^L|+\partial_z\kappa|\partial_t\eta^L|^2.
\end{eqnarray*}
Combining \eqref{pressure-estimate}, we have the desired formula.
\begin{remark}
We can rephrase the commutativity of $\partial_{R_1}$ and $\partial_z$ as 
\begin{equation*}
[\partial_z,\partial_{R_1}]=\partial_z\partial_{R_1}-\partial_{R_1}\partial_z=0,
\end{equation*}
where $[\cdot,\cdot]$ is the Lie braket (c.f. the Frobenius theorem, see Chapter 19 in \cite{L} for example).
For the previous studies using this property, see Chan-Czubak-Y \cite[Section 2.5]{CCY}
and Lichtenfelz-Y \cite{LY}, more originally, see Ma-Wang \cite[(3.7)]{MW}.
\end{remark}

\begin{remark}
Since $\nabla p\cdot b=\partial_{R_2} p\equiv 0$, then
\begin{equation*}
\begin{split}
-\partial_{R_2}(\nabla p\cdot \tau)
&=-\partial_{R_2}\partial_{z} p-T\partial_{R_1} p\\
(\text{commute}\ \partial_{R_2}\ \text{and}\  \partial_{z})&=
-T(\nabla p\cdot n)
=
T\kappa|\partial_t\eta^L|^2
\end{split}
\end{equation*}
for $R_1=R_2=0$.
However this formula is useless, since it coincides with \eqref{b}.
\end{remark}

\section{Proof of Corollary \ref{typical-strain-tensor}}
For any $x\in\{x: |x_2|>|x_1|\}\cap \{x_1\not=0\}$, let us set 
\begin{equation*}
\eta^L(t,x)=
\begin{pmatrix}
re^{t+t_0}\\
re^{-(t+t_0)}\\
x_{3}
\end{pmatrix},
\end{equation*}
where $r:=\sqrt{x_1x_2}$ and $t_0=\frac{1}{2}\log(x_1/x_2)$.
In this case  we see
\begin{equation*}
\partial_zt:=|\partial_t\eta^L|^{-1}=\frac{1}{r(e^{2t}+e^{-2t})^{1/2}}=\frac{1}{r\sqrt 2(\cosh 2t)^{1/2}}
\end{equation*}
and 
\begin{equation*}
\partial_z^2t=-\frac{(\sinh 2t)\partial_zt}{r\sqrt 2(\cosh 2t)^{3/2}}
=-\frac{\sinh 2t}{2r^2(\cosh 2t)^2} 
=-\frac{\tanh 2t}{2r^2\cosh 2t}.
\end{equation*}
On the other hand,
\begin{equation*}
\begin{split}
\kappa n&=\partial_t^2\eta^L(\partial_zt)^2+\partial_t\eta^L \partial_z^2t\\
&=
\frac{1}{2r^2\cosh 2t}
\begin{pmatrix}
re^{t}\\
re^{-t}\\
0
\end{pmatrix}
-\frac{\tanh 2t}{2r^2\cosh 2t}
\begin{pmatrix}
re^{t}\\
-re^{-t}\\
0
\end{pmatrix}.
\end{split}
\end{equation*}

Thus 
\begin{equation*}
\begin{split}
\kappa^2
&=
\frac{1}{2r^2\cosh 2t}-\frac{(\tanh 2t)^2}{r^2\cosh 2t}+\frac{(\tanh 2t)^2}{2r^2\cosh 2t}
=\frac{1}{2r^2(\cosh 2t)^3}
\end{split}
\end{equation*}
and then
\begin{equation*}
\begin{split}
\kappa&=\frac{1}{\sqrt 2r(\cosh 2t)^{3/2}},\quad
\partial_z\kappa=
-\frac{3\tanh 2t}{2r^2(\cosh 2t)^2},
\end{split}
\end{equation*}
\begin{equation*}
\begin{split}
|\partial_t\eta^L|^2=2r^2\cosh 2t\quad\text{and}\quad
\partial_t|\partial_t\eta^L|
=\sqrt 2r(\tanh 2t)^{1/2}(\sinh 2t)^{1/2}.
\end{split}
\end{equation*}
Combining the above calculations, we have the following desired formula:
\begin{equation*}
\begin{split}
\partial_z(\kappa|\partial_t\eta^L|^2)=
2\kappa \partial_t|\partial_t\eta^L|+\partial_z\kappa|\partial_t\eta^L|^2
=
-\frac{\tanh 2t}{\cosh 2t}.
\end{split}
\end{equation*}

\vspace{0.5cm}
\noindent
{\bf Acknowledgments.}\ 
Research of YS was partially supported by Grant-in-Aid for JSPS Fellows, Japan Society
for the Promotion of Science (JSPS).
Research of  TY  was partly supported by the JSPS Grants-in-Aid for Scientific
Research  20H01819.
\bibliographystyle{amsplain} 

\begin{thebibliography}{10} 
















\bibitem{CCY}
C-H. Chan, M. Czubak and T. Yoneda, 
\emph{An ODE for boundary layer separation on a sphere and a hyperbolic space,} 
 Physica D, \textbf{282} (2014) 34-38.

































\bibitem{Goto-2008}
S.~Goto,
\textit{A physical mechanism of the energy cascade in homogeneous isotropic
  turbulence},
J. Fluid Mech. \textbf{605} (2008) 355--366.




\bibitem{GSK}
S. Goto, Y. Saito, and G. Kawahara,
\textit{Hierarchy of antiparallel vortex tubes in spatially periodic turbulence at high Reynolds numbers},
Phys. Rev. Fluids \textbf{2} (2017) 064603. 

\bibitem{HSD}
P. E. Hamlington, J. Schumacher and W. J. A. Dahm,
\textit{Direct assessment of vorticity
alignment with local and nonlocal strain rates in turbulent flows}, Phys. Fluids \textbf{20} (2008) 111703.


\bibitem{JY1}
I.-J. Jeong and T. Yoneda, 
\textit{Enstrophy dissipation and vortex thinning for the incompressible 2D Navier-Stokes equations}, Nonlinearity \textbf{34} (2021) 1837.

\bibitem{JY2}
I.-J. Jeong and T. Yoneda, 
\textit{Vortex stretching and enhanced dissipation for the incompressible 3D Navier-Stokes equations}, Math. Annal. \textbf{380} (2021) 2041-2072.

\bibitem{JY3}
I.-J. Jeong and T. Yoneda, 
\textit{Quasi-streamwise vortices and enhanced dissipation for the incompressible 3D Navier-Stokes equations},  Proceedings of AMS \textbf{150} (2022) 1279-1286.














\bibitem{L}
J. M. Lee,
\emph{Introduction to smooth manifolds (second edition)}, 
Graduate Texts in Mathematics 218,
Springer, 2012.



\bibitem{LY}
L. Lichtenfelz and T. Yoneda, A local instability mechanism of the Navier-Stokes flow with swirl on the no-slip flat boundary, J. Math. Fluid Mech. \textbf{21} (2019) 20.




\bibitem{MW}
T. Ma and S. Wang,
\emph{
Boundary layer separation and structural bifurcation for 2-D incompressible fluid flows.
 Partial differential equations and applications,} Discrete Contin. Dyn. Syst. \textbf{10} (2004)  459--472.




\bibitem{MP}
 G. Misio{\l}ek and S. Preston, 
\textit{Fredholm properties of Riemannian exponential maps on diffeomorphism groups}, Invent. math. \textbf{179} (2010) 191--227.







\bibitem{Motoori-2019}
Y.~Motoori and S.~Goto, 
\textit{Generation mechanism of a hierarchy of vortices in a turbulent
  boundary layer}, 
J. Fluid Mech. \textbf{865} (2019) 1085--1109.

\bibitem{Motoori-2021}
Y.~Motoori and S.~Goto, 
\textit{Hierarchy of coherent structures and real-space energy transfer in
  turbulent channel flow}, 
 J. Fluid Mech. \textbf{911} (2021) A27.

















\bibitem{TGOY}
T. Tsuruhashi, S. Goto,  S. Oka and  T. Yoneda, 
\textit{Self-similar hierarchy of coherent tubular vortices in turbulence}, to appear in Philosophical Transactions A.


\bibitem{YGT}
T. Yoneda, S. Goto and T. Tsuruhashi, 
\textit{Mathematical reformulation of the Kolmogorov-Richardson energy cascade in terms of vortex stretching}, Nonlinearity \textbf{34} (2021) 1837.








\end{thebibliography}

\end{document}